# Basic Properties of Strong Mixing Conditions. A Survey and Some Open Questions[*]

**Richard C. Bradley**

*e-mail:* `bradleyr@indiana.edu`

**Abstract:** This is an update of, and a supplement to, a 1986 survey paper by the author on basic properties of strong mixing conditions.

This is an update of, and a supplement to, the author's earlier survey paper [18] on basic properties of strong mixing conditions. That paper appeared in 1986 in a book containing survey papers on various types of dependence conditions and the limit theory under them. The survey here will include part (but not all) of the material in [18], and will also describe some relevant material that was not in that paper, especially some new discoveries and developments that have occurred since that paper was published. (Much of the new material described here involves "interlaced" strong mixing conditions, in which the index sets are not restricted to "past" and "future.") At various places in this survey, open problems will be posed.

There is a large literature on basic properties of strong mixing conditions. A survey such as this cannot do full justice to it. Here are a few references on important topics not covered in this survey. For the approximation of mixing sequences by martingale differences, see e.g. the book by Hall and Heyde [80]. For the direct approximation of mixing random variables by independent ones, see e.g. [43, Chapter 16], [62], [111], [131], [136, Chapter 5]. For some "coupling" properties connected with the "absolute regularity" condition, see [2, Chapter 4]. For mixing properties of linear processes, see [71] and [143]. For some very strong mixing properties of one-dimensional Gibbs states, see e.g. [63] and [66]. For some very strong mixing properties of a well known "continued fraction" process, see e.g. [100] and [132]. For a broad survey on the connections between dynamical systems and strong mixing conditions, see the survey paper by Denker [63].

This survey here is organized as follows:
1. Measures of dependence
2. Some strong mixing conditions

---

[*]This is an original survey paper.

3. Markov chains
4. General behavior of the dependence coefficients
5. Independent pairs of $\sigma$–fields
6. Linear dependence conditions
7. Gaussian sequences
8. Random fields

# 1. The measures of dependence

## 1.1. Definitions and some basic properties

In what follows, expressions such as $\sup_{q\in Q, s\in S} h(q,s)$ will often be written as $\sup h(q,s)$, $q \in Q$, $s \in S$.

Throughout this paper, the probability space is $(\Omega, \mathcal{F}, P)$.

For any $\sigma$-field $\mathcal{A} \subset \mathcal{F}$, let $\mathcal{L}^2_{\text{real}}(\mathcal{A})$ denote the space of (equivalence classes of) square-integrable, $\mathcal{A}$-measurable (real-valued) random variables.

For any two $\sigma$-fields $\mathcal{A}$ and $\mathcal{B} \subset \mathcal{F}$, define the following eight measures of dependence:

$$\alpha(\mathcal{A}, \mathcal{B}) := \sup |P(A \cap B) - P(A)P(B)|, \quad A \in \mathcal{A}, B \in \mathcal{B}; \tag{1.1}$$

$$\phi(\mathcal{A}, \mathcal{B}) := \sup |P(B|A) - P(B)|, \quad A \in \mathcal{A}, B \in \mathcal{B}, P(A) > 0; \tag{1.2}$$

$$\psi(\mathcal{A}, \mathcal{B}) := \sup \left| \frac{P(A \cap B)}{P(A)P(B)} - 1 \right|, A \in \mathcal{A}, B \in \mathcal{B}, P(A) > 0, P(B) > 0; \tag{1.3}$$

$$\rho(\mathcal{A}, \mathcal{B}) := \sup |\text{Corr}(f, g)|, \quad f \in \mathcal{L}^2_{\text{real}}(\mathcal{A}), g \in \mathcal{L}^2_{\text{real}}(\mathcal{B}); \tag{1.4}$$

$$\beta(\mathcal{A}, \mathcal{B}) := \sup \frac{1}{2} \sum_{i=1}^{I} \sum_{j=1}^{J} |P(A_i \cap B_j) - P(A_i)P(B_j)| \tag{1.5}$$

where the supremum is taken over all pairs of (finite) partitions $\{A_1, \ldots, A_I\}$ and $\{B_1, \ldots, B_J\}$ of $\Omega$ such that $A_i \in \mathcal{A}$ for each $i$ and $B_j \in \mathcal{B}$ for each $j$;

$$\psi^*(\mathcal{A}, \mathcal{B}) := \sup \frac{P(A \cap B)}{P(A)P(B)}; \quad A \in \mathcal{A}, B \in \mathcal{B}, P(A) > 0, P(B) > 0; \tag{1.6}$$

$$\psi'(\mathcal{A}, \mathcal{B}) := \inf \frac{P(A \cap B)}{P(A)P(B)}; \quad A \in \mathcal{A}, B \in \mathcal{B}, P(A) > 0, P(B) > 0; \tag{1.7}$$

$$I(\mathcal{A}, \mathcal{B}) := \sup \sum_{i=1}^{I} \sum_{j=1}^{J} P(A_i \cap B_j) \log \left( \frac{P(A_i \cap B_j)}{P(A_i)P(B_j)} \right) \tag{1.8}$$

where the supremum is taken over all pairs of (finite) partitions $\{A_1, \ldots, A_I\}$ and $\{B_1, \ldots, B_J\}$ of $\Omega$ such that $\mathcal{A}_i \in \mathcal{A}$ for each $i$ and $B_j \in \mathcal{B}$ for each $j$. In (1.8) and in what follows, $0/0 := 0$ and $0 \log 0 := 0$.

The "maximal correlation" coefficient $\rho(\mathcal{A}, \mathcal{B})$ was first studied in the papers [74], [77], [86], [87], [109], in statistical contexts that had no particular connection with stochastic processes. The "coefficient of information" $I(\mathcal{A}, \mathcal{B})$ (along with the related notion of entropy) was developed in papers such as [78] and [146].

The following inequalities give the ranges of possible values (including $\infty$ in some cases) of those measures of dependence:

$$
\begin{aligned}
&0 \leq \alpha(\mathcal{A}, \mathcal{B}) \leq 1/4, \quad 0 \leq \phi(\mathcal{A}, \mathcal{B}) \leq 1, \quad 0 \leq \psi(\mathcal{A}, \mathcal{B}) \leq \infty, \\
&0 \leq \rho(\mathcal{A}, \mathcal{B}) \leq 1, \quad 0 \leq \beta(\mathcal{A}, \mathcal{B}) \leq 1, \quad 1 \leq \psi^*(\mathcal{A}, \mathcal{B}) \leq \infty, \\
&0 \leq \psi'(\mathcal{A}, \mathcal{B}) \leq 1, \quad 0 \leq I(\mathcal{A}, \mathcal{B}) \leq \infty.
\end{aligned} \quad (1.9)
$$

Each of the following equalities is equivalent to the condition that $\mathcal{A}$ and $\mathcal{B}$ are independent:

$$
\begin{aligned}
&\alpha(\mathcal{A}, \mathcal{B}) = 0, \quad \phi(\mathcal{A}, \mathcal{B}) = 0, \quad \psi(\mathcal{A}, \mathcal{B}) = 0, \quad \rho(\mathcal{A}, \mathcal{B}) = 0, \quad \beta(\mathcal{A}, \mathcal{B}) = 0, \\
&\psi^*(\mathcal{A}, \mathcal{B}) = 1, \quad \psi'(\mathcal{A}, \mathcal{B}) = 1, \quad I(\mathcal{A}, \mathcal{B}) = 0.
\end{aligned} \quad (1.10)
$$

These measures of dependence satisfy the following inequalities:

$$
\begin{aligned}
2\alpha(\mathcal{A}, \mathcal{B}) &\leq \beta(\mathcal{A}, \mathcal{B}) \leq \phi(\mathcal{A}, \mathcal{B}) \leq (1/2)\psi(\mathcal{A}, \mathcal{B}); & (1.11) \\
4\alpha(\mathcal{A}, \mathcal{B}) &\leq \rho(\mathcal{A}, \mathcal{B}) \leq \psi(\mathcal{A}, \mathcal{B}); & (1.12) \\
\rho(\mathcal{A}, \mathcal{B}) &\leq 2[\phi(\mathcal{A}, \mathcal{B})]^{1/2}[\phi(\mathcal{B}, \mathcal{A})]^{1/2} \leq 2[\phi(\mathcal{A}, \mathcal{B})]^{1/2}; & (1.13) \\
\phi(\mathcal{A}, \mathcal{B}) &\leq 1 - 1/\psi^*(\mathcal{A}, \mathcal{B}) \leq \psi^*(\mathcal{A}, \mathcal{B}) - 1; & (1.14) \\
\phi(\mathcal{A}, \mathcal{B}) &\leq 1 - \psi'(\mathcal{A}, \mathcal{B}); & (1.15) \\
\psi(\mathcal{A}, \mathcal{B}) &= \max\{\psi^*(\mathcal{A}, \mathcal{B}) - 1, 1 - \psi'(\mathcal{A}, \mathcal{B})\}; & (1.16) \\
I(\mathcal{A}, \mathcal{B}) &\leq \psi^*(\mathcal{A}, \mathcal{B}) \log \psi^*(\mathcal{A}, \mathcal{B}); & (1.17) \\
\beta(\mathcal{A}, \mathcal{B}) &\leq [I(\mathcal{A}, \mathcal{B})]^{1/2}. & (1.18)
\end{aligned}
$$

The first inequality in (1.13) was shown by Peligrad [124] with an extension of the arguments used by Cogburn [57] and Ibragimov [91] to show the inequality $\rho(\mathcal{A}, \mathcal{B}) \leq 2[\phi(\mathcal{A}, \mathcal{B})]^{1/2}$ (see also Doob [70, p. 222, Lemma 7.1]). Denker and Keller [65] independently proved the similar inequality $\rho(\mathcal{A}, \mathcal{B}) \leq 2 \cdot \max\{\phi(\mathcal{A}, \mathcal{B}), \phi(\mathcal{B}, \mathcal{A})\}$. Equation (1.18) essentially goes back to Volkonskii and Rozanov [162], [163]. The other inequalities are elementary.

### *1.2. A brief look at norms*

By [164, Theorem 1.1] and a simple argument, for any two $\sigma$-fields $\mathcal{A}$ and $\mathcal{B}$,

$$
\rho(\mathcal{A}, \mathcal{B}) = \sup \frac{|Efg - EfEg|}{\|f\|_2 \|g\|_2} \quad (1.19)
$$

where the supremum is taken over all pairs of (not necessarily centered) *complex-valued* absolutely square-integrable random variables $f$ and $g$ such that $f$ is $\mathcal{A}$-measurable and $g$ is $\mathcal{B}$-measurable.

In a similar spirit (we omit the details here), one can express a broad class of measures of dependence, including some of the ones in section 1.1, as "norms" (of various kinds, with various parameters) of the bilinear form "covariance." Then one can apply results and techniques from functional analysis, including interpolation theory (see [1] or [8]), to efficiently compare those measures of dependence and derive a broad class of "covariance inequalities." See e.g. [138], [106], [45], [46], and [47]. The last two papers also gave some constructions to help establish "equivalence classes" of those measures of dependence.

Here we shall just look at one open problem arising from the last paper.

Suppose $B$ is a (say) real Banach space, with norm $\|.\|_B$. Let $B^*$ denote its dual space (the space of real bounded linear functionals on $B$), with its usual norm, denoted $\|\cdot\|_{B^*}$. For any $x \in B$ and any $y \in B^*$, the real number $y(x)$ will be denoted $\langle x, y \rangle$ (by analogy with the case of a Hilbert space, where a bounded linear functional is an inner product with a fixed element).

Following (essentially) the notations in [47], for any two $\sigma$-fields $\mathcal{A}$ and $\mathcal{B}$, define the measure of dependence

$$R^B_{\infty,\infty}(\mathcal{A},\mathcal{B}) := \sup \frac{|E\langle X, Y\rangle - \langle EX, EY\rangle|}{\|\,\|X\|_B\,\|_\infty \cdot \|\,\|Y\|_{B^*}\,\|_\infty}$$

where the supremum is taken over all pairs of simple random variables $X$ and $Y$ such that $X$ is $B$-valued and $\mathcal{A}$-measurable and $Y$ is $B^*$-valued and $\mathcal{B}$-measurable.

Let us say that two measures of dependence are "equivalent" if each one becomes arbitrarily small as the other becomes sufficiently small. By Dehling and Philipp [61, Lemma 2.2], for any nontrivial Hilbert space $H$, the measure of dependence $R^H_{\infty,\infty}(\cdot,\cdot)$ is equivalent to $\alpha(\cdot,\cdot)$. With an elementary construction, it was shown in [47, Theorem 3.1] that for $B = \ell^\infty$, $\ell^1$ or $c_0$ (the subspace of $\ell^\infty$ consisting of sequences that converge to 0), $R^B_{\infty,\infty}(\mathcal{A},\mathcal{B}) = 2\beta(\mathcal{A},\mathcal{B})$.

The following question remains open:

**Question 1**: *For an arbitrary nontrivial real Banach space $B$, is the measure of dependence $R^B_{\infty,\infty}(\cdot,\cdot)$ equivalent to one of the two measures of dependence $\alpha(\cdot,\cdot)$, $\beta(\cdot,\cdot)$?*

The author [37] showed that for any given $c \in (1,\infty)$, if one restricts the notion of "equivalence" to pairs of $\sigma$–fields $\mathcal{A}, \mathcal{B}$ such that $\psi^*(\mathcal{A},\mathcal{B}) \leq c$, then the answer to this question is affirmative (even for a class of measures of dependence much broader than just the ones $R^B_{\infty,\infty}(\cdot,\cdot)$).

## 2. Strong mixing conditions

### 2.1. Strong mixing conditions based on "past" and "future"

Suppose $X := (X_k, k \in \mathbb{Z})$ is a (not necessarily stationary) sequence of random variables. For $-\infty \leq J \leq L \leq \infty$, define the $\sigma$-field

$$\mathcal{F}^L_J := \sigma(X_k, J \leq k \leq L \ (k \in \mathbb{Z})). \tag{2.1}$$

Here and below, the notation $\sigma(\ldots)$ means the $\sigma$-field $\subset \mathcal{F}$ generated by $(\ldots)$. For each $n \geq 1$, define the following dependence coefficients:

$$\begin{aligned}
\alpha(n) &:= \sup_{j\in\mathbb{Z}} \alpha(\mathcal{F}_{-\infty}^{j}, \mathcal{F}_{j+n}^{\infty}); \\
\phi(n) &:= \sup_{j\in\mathbb{Z}} \phi(\mathcal{F}_{-\infty}^{j}, \mathcal{F}_{j+n}^{\infty}); \\
\psi(n) &:= \sup_{j\in\mathbb{Z}} \psi(\mathcal{F}_{-\infty}^{j}, \mathcal{F}_{j+n}^{\infty}); \\
\rho(n) &:= \sup_{j\in\mathbb{Z}} \rho(\mathcal{F}_{-\infty}^{j}, \mathcal{F}_{j+n}^{\infty}); \\
\beta(n) &:= \sup_{j\in\mathbb{Z}} \beta(\mathcal{F}_{-\infty}^{j}, \mathcal{F}_{j+n}^{\infty}); \\
\psi^{*}(n) &:= \sup_{j\in\mathbb{Z}} \psi^{*}(\mathcal{F}_{-\infty}^{j}, \mathcal{F}_{j+n}^{\infty}); \\
\psi'(n) &:= \inf_{j\in\mathbb{Z}} \psi'(\mathcal{F}_{-\infty}^{j}, \mathcal{F}_{j+n}^{\infty}); \quad \text{and} \\
I(n) &:= \sup_{j\in\mathbb{Z}} I(\mathcal{F}_{-\infty}^{j}, \mathcal{F}_{j+n}^{\infty}).
\end{aligned} \qquad (2.2)$$

(Note the "inf" in the definition of $\psi'(n)$.)

The random sequence $X$ is said to be
"strongly mixing" (or "$\alpha$-mixing") if $\alpha(n) \to 0$ as $n \to \infty$,
"$\phi$-mixing" if $\phi(n) \to 0$ as $n \to \infty$,
"$\psi$-mixing" if $\psi(n) \to 0$ as $n \to \infty$,
"$\rho$-mixing" if $\rho(n) \to 0$ as $n \to \infty$,
"absolutely regular" (or "$\beta$-mixing") if $\beta(n) \to 0$ as $n \to \infty$,
"$\psi^{*}$-mixing" if $\psi^{*}(n) \to 1$ as $n \to \infty$,
"$\psi'$-mixing" if $\psi'(n) \to 1$ as $n \to \infty$,
"information regular" if $I(n) \to 0$ as $n \to \infty$.

The strong mixing ($\alpha$-mixing) condition was introduced by Rosenblatt [137]. The $\phi$-mixing condition was introduced by Ibragimov [88], and was also studied by Cogburn [57]. The $\psi$-mixing condition had its origin in a paper by Blum, Hanson, and Koopmans [9] studying a different condition ("*-mixing") based on the same measure of dependence, and it took its present form in the paper of Philipp [130]. The $\rho$-mixing condition was introduced by Kolmogorov and Rozanov [104]. The absolute regularity condition was introduced by Volkonskii and Rozanov [162], [163], and was attributed there to Kolmogorov. The $\psi^{*}$-mixing and $\psi'$-mixing conditions are obvious "halves" of the $\psi$-mixing condition; their origins are hard to trace. The information regularity condition was introduced by Volkonskii and Rozanov [162], [163], and was (at least in spirit) attributed there to M.S. Pinsker.

In the special case where the sequence $X$ is strictly stationary, one has simply

$$\alpha(n) = \alpha(\mathcal{F}_{-\infty}^{0}, \mathcal{F}_{n}^{\infty}), \qquad (2.3)$$

and the same holds for the other dependence coefficients in (2.2).

For strictly stationary sequences $X$, the following is another (equivalent) formulation of the absolute regularity condition:

$$\forall \, \varepsilon > 0, \quad \exists \, n \geq 1, \quad \exists \, D \in \mathcal{F}^0_{-\infty} \quad \text{with } P(D) \geq 1 - \varepsilon, \text{ such that}$$
$$\forall \, A \in \mathcal{F}^0_{-\infty} \quad \text{such that } A \subset D \text{ and } P(A) > 0, \, \forall \, B \in \mathcal{F}^\infty_n,$$
$$\text{one has that } |P(B|A) - P(B)| \leq \varepsilon. \tag{2.4}$$

For strictly stationary, finite-state sequences, this formulation (in essence) was given by Friedman and Ornstein [76] under the name "weak Bernoulli condition."

**A caution on the terminology.** It needs to be kept in mind that two barely different phrases are used with quite different meanings: The phrase "strong mixing condition" (singular), or simply "strong mixing," refers to $\alpha$–mixing ($\alpha(n) \to 0$) as above. In contrast, the phrase "strong mixing conditions" (plural) refers to all mixing conditions that are at least as strong as (i.e. that imply) $\alpha$–mixing.

The latter phrase "strong mixing conditions" is intended to distinguish from a broad class of "mixing conditions" from ergodic theory that are weaker than $\alpha$–mixing. (See e.g. Petersen [129].)

From eqs. (1.11)—(1.18), one has the following "hierarchy" of these strong mixing conditions:
(a) $m$-dependence implies $\psi$-mixing.
(b) $\psi$-mixing implies $\psi^*$-mixing.
(c) $\psi$-mixing implies $\psi'$-mixing.
(d) $\psi^*$-mixing implies $\phi$-mixing.
(e) $\psi^*$-mixing implies information regularity.
(f) $\psi'$-mixing implies $\phi$-mixing.
(g) Information regularity implies absolute regularity.
(h) $\phi$-mixing implies absolute regularity.
(i) $\phi$-mixing implies $\rho$-mixing.
(j) Absolute regularity implies strong mixing.
(k) $\rho$-mixing implies strong mixing.

Aside from "transitivity," there are no other implications between these mixing conditions. (For more on that, including references to examples, see e.g. [42, Remark 5.23].)

### 2.2. The asymmetry of $\phi$-mixing

For any two to $\sigma$-fields $\mathcal{A}$ and $\mathcal{B}$, obviously $\alpha(\mathcal{A}, \mathcal{B}) = \alpha(\mathcal{B}, \mathcal{A})$. The same symmetry holds for the other measures of dependence in (1.1)–(1.8) except $\phi(\cdot, \cdot)$. If $\phi(\mathcal{A}, \mathcal{B})$ is "small," that does not imply that $\phi(\mathcal{B}, \mathcal{A})$ is "small."

Suppose $X := (X_k, k \in \mathbb{Z})$ is a strictly stationary random sequence. Then for each $n \geq 1$, $\phi(n) = \phi(\mathcal{F}^0_{-\infty}, \mathcal{F}^\infty_n)$. For each $n \geq 1$, define

$$\phi_{\text{rev}}(n) := \phi\left(\mathcal{F}^\infty_n, \mathcal{F}^0_{-\infty}\right). \tag{2.5}$$

(The subscript "rev" is an abbreviation of "reversed.") The sequence $X$ is said to be "time–reversed $\phi$-mixing" if $\phi_{\text{rev}}(n) \to 0$ as $n \to \infty$.

Rosenblatt [138, pp. 213–214] constructed some strictly stationary Markov chains that are $\phi$-mixing but not "time–reversed $\phi$-mixing."

## 2.3. Interlaced strong mixing conditions

Suppose $X := (X_k, k \in \mathbb{Z})$ is a (not necessarily stationary) random sequence. For each $n \geq 1$, define the following dependence coefficients:

$$\begin{align}
\alpha^*(n) &:= \sup \alpha(\sigma(X_k, k \in S), \sigma(X_k, k \in T)), \tag{2.6} \\
\rho^*(n) &:= \sup \rho(\sigma(X_k, k \in S), \sigma(X_k, k \in T)), \quad \text{and} \tag{2.7} \\
\beta^*(n) &:= \sup \beta(\sigma(X_k, k \in S), \sigma(X_k, k \in T)), \tag{2.8}
\end{align}$$

where in each of (2.6), (2.7), and (2.8), the supremum is taken over all pairs of nonempty, disjoint sets $S$ and $T \subset \mathbb{Z}$ such that

$$\text{dist}(S, T) := \min_{s \in S, t \in T} |s - t| \geq n. \tag{2.9}$$

In (2.9), it is understood that the two sets $S$ and $T$ can be "interlaced," i.e. with each set having elements between ones in the other set.

**Theorem 2.1** *Suppose $X := (X_k, k \in \mathbb{Z})$ is a strictly stationary random sequence.*
*(1) $\beta^*(n) \to 0$ as $n \to \infty$ if and only if $X$ is $m$–dependent.*
*(2) $\alpha^*(n) \to 0$ as $n \to \infty$ if and only if $\rho^*(n) \to 0$ as $n \to \infty$.*

Part (1) was shown in [22, Theorem 1 and Remarks 2 and 3], as part of a broader result for strictly stationary random fields. The main insight behind it came from examples that were presented by Dobrushin [68, p. 205] and Zhurbenko [166, p. 8] to show that for strictly stationary random fields, a seemingly natural formulation of a "$\phi$-mixing" condition turns out to be much stronger than it appears at first sight.

Part (2) was shown in [26, Theorem 1 and Remarks 1, 2, and 3], as part of a broader result for strictly stationary random fields. Versions of (2) for strictly stationary random fields had also been noticed in the 1980's by probabilists (faculty and students) at Moscow State University (in Moscow, Russia), but apparently they never published anything on that. (According to Zhurbenko [167], there may have been some uncertainty there about the statements or proofs.) For stationary Gaussian random fields, a version of (2) had been proved by Rosenblatt [140], with an adaptation of an argument of Kolmogorov and Rozanov [104].

Primarily because of Theorem 2.1, only one of the three dependence coefficients in (2.6)–(2.8) will be used in a formal definition of a mixing condition:

A given random sequence $X := (X_k, k \in \mathbb{Z})$ is said to be "$\rho^*$-mixing" if $\rho^*(n) \to 0$ as $n \to \infty$.

Under strict stationarity, the $\rho^*$-mixing condition goes back at least to Stein [154] and (in the equivalent form $\alpha^*(n) \to 0$) to Rosenblatt [140]. (The latter paper involved random fields.)

For (say) strictly stationary sequences, the exact location of $\rho^*$-mixing in the "hierarchy" at the end of section 2.1 has not yet been established. Obviously $\rho^*$-mixing implies $\rho$-mixing. From known examples (see e.g. [34, Example 6.4] or [42, Example 7.16]), one has that $\rho^*$-mixing does not imply absolute regularity. Also, $m$-dependence implies $\rho^*$-mixing (this is elementary); and $\rho$-mixing does not imply $\rho^*$-mixing (see e.g. [36]). But that does not give the whole picture. In particular, the following question remains unsolved:

**Question 2**: *If a given random sequence is $\phi$-mixing, does it follow that it is $\rho^*$-mixing?*

An affirmative answer to this question was conjectured by the author [31, p. 226]. Walter Philipp [133] said he thinks the answer is negative instead.

In the 1960's, I.A. Ibragimov conjectured that if a given strictly stationary sequence $X := (X_k, k \in \mathbb{Z})$ is $\phi$-mixing, has finite second moments, and satisfies $\text{Var}(X_1 + \ldots + X_n) \to \infty$ as $n \to \infty$, then it satisfies a CLT. (See [95, p. 393, problem (3)].) Iosifescu [99, p. 56] conjectured that under the same hypothesis, a weak invariance principle holds. These conjectures remain unsolved. Peligrad [125, Corollary 2.2 and p. 1305, lines 7–8] has confirmed them under the augmented hypothesis $\liminf_{n\to\infty} n^{-1}\text{Var}(X_1 + \ldots + X_n) > 0$. If $\phi$-mixing implies $\rho^*$-mixing, then these conjectures of Ibragimov and Iosifescu (under just $\text{Var}(X_1 + \ldots + X_n) \to \infty$) would follow immediately from known corresponding results under $\rho^*$-mixing (see e.g. [24] and [125]).

### 2.4. $\rho$-mixing except on small sets, and a "two–part" mixing condition

For any event $D$ with $P(D) > 0$, and any two $\sigma$-fields $\mathcal{A}$ and $\mathcal{B} \subset \mathcal{F}$, let $\rho_{P(.|D)}(\mathcal{A}, \mathcal{B})$ denote the maximal correlation coefficient between $\mathcal{A}$ and $\mathcal{B}$ with respect to the (conditional) probability measure on $(\Omega, \mathcal{F})$ given by $F \mapsto P(F|D)$, $F \in \mathcal{F}$. For any two $\sigma$-fields $\mathcal{A}$ and $\mathcal{B} \subset \mathcal{F}$, define the measure of dependence

$$\rho_{\text{cond}}(\mathcal{A}, \mathcal{B}) := \inf\{\varepsilon \in (0,1) : \exists D \in \mathcal{A} \text{ such that}$$
$$P(D) \geq 1 - \varepsilon \text{ and } \rho_{P(.|D)}(\mathcal{A}, \mathcal{B}) \leq \varepsilon\}. \quad (2.10)$$

(If no such $\varepsilon \in (0,1)$ exists, let $\rho_{\text{cond}}(\mathcal{A}, \mathcal{B}) := 1$.) The subscript "cond" is an abbreviation for "conditional." Equation (2.10) is an analog of the formulation of the absolute regularity condition given in eq. (2.4).

In [23, Proposition 2.1 and subsequent comments, and p. 219, lines 17–19] (see [43, Note 1 after Lemma 24.11]) it was shown that (i) $\alpha(\mathcal{A}, \mathcal{B}) \leq 4\rho_{\text{cond}}(\mathcal{A}, \mathcal{B})$ for any two $\sigma$-fields $\mathcal{A}$ and $\mathcal{B}$, and (ii) there exists $\varepsilon_0 > 0$ such that for any two given $\sigma$-fields $\mathcal{A}$ and $\mathcal{B}$, if $\alpha(\mathcal{A}, \mathcal{B}) \leq \varepsilon_0$ then $\rho_{\text{cond}}(\mathcal{A}, \mathcal{B}) \leq \alpha^{1/12}(\mathcal{A}, \mathcal{B})$. (The exponent 1/12 can be improved with trivial refinements of the arithmetic in that argument.) This gives an "equivalence" between the two measures of dependence $\alpha(.,.)$ and $\rho_{\text{cond}}(\mathcal{A}, \mathcal{B})$, in that each one becomes arbitrarily small as the other becomes sufficiently small. Theorem 2.2 below is an immediate corollary. By a simple argument, comments (i) and (ii) above and the resulting "equivalence" and Theorem 2.2 below, all hold if in (2.10) the condition $D \in \mathcal{A}$ is replaced by $D \in \mathcal{F}$.

**Theorem 2.2** *For a given (not necessarily stationary) random sequence $X := (X_k, k \in \mathbb{Z})$, the following two conditions are equivalent:*
*(a) $\alpha(n) \to 0$ as $n \to \infty$ (strong mixing).*
*(b) $\rho_{\text{cond}}(n) := \sup_{j \in \mathbb{Z}} \rho_{\text{cond}}(\mathcal{F}_{-\infty}^j, \mathcal{F}_{j+n}^\infty) \to 0$ as $n \to \infty$.*

For strictly stationary sequences, this theorem was formulated in [23, Theorem 1.3].

Condition (b) in Theorem 2.2 might be referred to as a condition of "$\rho$-mixing except on small sets." Such notions of "$\rho$-mixing except on small sets," and related notions of a "two-part mixing condition" (a "hybrid" of strong mixing and $\rho$-mixing), both arose in a conversation between Magda Peligrad and Enrico Presutti in the spring of 1983, and they were treated in the paper by Peligrad and the author [48].

The "two-part mixing condition" is as follows, formulated here for a given strictly stationary sequence $X := (X_k, k \in \mathbb{Z})$: There exist nonincreasing sequences $(a_1, a_2, a_3, \ldots)$ and $(z_1, z_2, z_3, \ldots)$ of numbers in $[0, 1]$ with $a_n \to 0$ and $z_n \to 0$ as $n \to \infty$, such that

$$\forall\, n \geq 1, \quad \forall\, A \in \mathcal{F}_{-\infty}^0, \quad \forall\, B \in \mathcal{F}_n^\infty,$$
$$|P(A \cap B) - P(A)P(B)| \leq a_n + z_n[P(A)P(B)]^{1/2}. \qquad (2.11)$$

This condition is equivalent to strong mixing: If $X$ is strongly mixing, then (2.11) holds for $a_n := \alpha(n)$ and $z_n := 0$; and conversely, (2.11) implies strong mixing with $\alpha(n) \leq a_n + z_n$. However, in general the "mixing rates" do not "match." For example, by [48, Theorem 3], for a given (large) $\theta > 0$ and a given sequence $(z_1, z_2, z_3, \ldots)$ of positive numbers converging to 0 (very slowly), there exists a strictly stationary sequence $X := (X_k, k \in \mathbb{Z})$ such that (2.11) holds with $a_n \asymp n^{-\theta}$ as $n \to \infty$, and $\alpha(n) \asymp z_n$ as $n \to \infty$ and $\rho(n) = 1$ for all $n \geq 1$. Consequently, central limit theorems under (2.11) (such as the very sharp CLT in Peligrad [127]) cannot be derived from the usual CLT's under strong mixing (such as in [91], [72], [110] or [126]) or under $\rho$-mixing (such as in [94] or [126]). Bryc and Peligrad [52] showed that if Tukey's [157] "3R" (or "running median") smoother is applied to a strictly stationary, $\rho$-mixing sequence $X := (X_k, k \in \mathbb{Z})$, then the (strictly stationary) "output" sequence satisfies (2.11) with $a_n \to 0$ at least exponentially fast (and $z_n \to 0$, possibly very slowly). Under reasonable moment conditions, Bryc and Peligrad then applied the CLT in [127] under (2.11) to that "output" sequence.

### 2.5. Tail $\sigma$-fields

Suppose $X := (X_k, k \in \mathbb{Z})$ is a (not necessarily stationary) random sequence. For this sequence $X$,
the "past tail $\sigma$-field" is $\mathcal{T}_{\text{past}} := \bigcap_{n=1}^{\infty} \mathcal{F}_{-\infty}^{-n}$,
the "future tail $\sigma$-field" is $\mathcal{T}_{\text{future}} := \bigcap_{n=1}^{\infty} \mathcal{F}_n^{\infty}$,
the "double tail $\sigma$-field" is $\mathcal{T}_{\text{double}} := \bigcap_{n=1}^{\infty} (\mathcal{F}_{-\infty}^{-n} \vee \mathcal{F}_n^{\infty})$.

Obviously $\mathcal{T}_{\text{past}} \subset \mathcal{T}_{\text{double}}$ and $\mathcal{T}_{\text{future}} \subset \mathcal{T}_{\text{double}}$.

A $\sigma$-field $\mathcal{A} \subset \mathcal{F}$ is said to be "trivial" if $P(A) = 0$ or 1 for every $A \in \mathcal{A}$.

By elementary arguments, the following implications hold:
(a) If $\alpha(n) \to 0$ as $n \to \infty$, then $\mathcal{T}_{\text{past}}$ and $\mathcal{T}_{\text{future}}$ are each trivial.
(b) If $\rho(n) < 1$ for some $n \geq 1$, then $\mathcal{T}_{\text{past}}$ and $\mathcal{T}_{\text{future}}$ are trivial.
(c) If $\psi^*(n) < 2$ for some $n \geq 1$, then $\mathcal{T}_{\text{past}}$ and $\mathcal{T}_{\text{future}}$ are trivial.
(d) If $\psi'(n) > 0$ for some $n \geq 1$, then $\mathcal{T}_{\text{past}}$ and $\mathcal{T}_{\text{future}}$ are trivial.
(e) If $\beta(n) \to 0$ as $n \to \infty$, then $\mathcal{T}_{\text{double}}$ is trivial.
(f) If $\rho^*(n) < 1$ for some $n \geq 1$, then $\mathcal{T}_{\text{double}}$ is trivial.

Even in the case where $X$ is strictly stationary, either $\mathcal{T}_{\text{past}}$ or $\mathcal{T}_{\text{future}}$ can be trivial without the other being trivial. For an old, classic, well known example, a simple autoregressive process of order 1, see e.g. [141, p. 267], [18, Example 6.2] or [42, Example 2.15].

A given strictly stationary random sequence $X := (X_k, k \in \mathbb{Z})$ is said to be "mixing (in the ergodic–theoretic sense)" if

$$\forall\, A, B \in \mathcal{R}^{\mathbb{Z}}, \quad P(X \in A \cap T^{-n} B) \longrightarrow P(X \in A) P(X \in B) \text{ as } n \to \infty. \quad (2.12)$$

Here $T$ is the usual shift operator on $\mathbb{R}^{\mathbb{Z}}$. That is, for $\omega := (\omega_k, k \in \mathbb{Z}) \in \mathbb{R}^{\mathbb{Z}}$, the element $T\omega \in \mathbb{R}^{\mathbb{Z}}$ is given by $(T\omega)_k = \omega_{k+1}$, $k \in \mathbb{Z}$. Also, here and below, $\mathcal{R}$ (resp. $\mathcal{R}^{\mathbb{Z}}$) denotes the Borel $\sigma$-field on $\mathbb{R}$ (resp. $\mathbb{R}^{\mathbb{Z}}$).

By a well known result of Vinokurov [161] (see e.g. [95, Theorem 17.1.1] or [42, Theorem 2.14]), if $X := (X_k, k \in \mathbb{Z})$ is a strictly stationary sequence such that either $\mathcal{T}_{\text{past}}$ or $\mathcal{T}_{\text{future}}$ is trivial, then $X$ is mixing (in the ergodic–theoretic sense).

As is also well known, if a given strictly stationary sequence is mixing (in the ergodic–theoretic sense), then it is ergodic.

A treatment of other related concepts (such as "weak mixing") in ergodic theory can be found e.g. in Petersen [129].

### 2.6. Bilaterally deterministic sequences

A random sequence $X := (X_k, k \in \mathbb{Z})$ (stationary or not) is said to be "bilaterally deterministic" if $\mathcal{T}_{\text{double}} \stackrel{.}{=} \mathcal{F}_{-\infty}^{\infty}$ —that is, modulo null–sets, the double tail $\sigma$-field gives the history of the entire sequence.

Olshen [118] constructed a strictly stationary sequence $X := (X_k, k \in \mathbb{Z})$ such that $\mathcal{T}_{\text{past}}$ and $\mathcal{T}_{\text{future}}$ are each trivial but $\mathcal{T}_{\text{double}}$ is not trivial. (In that construction $\mathcal{T}_{\text{double}}$ was not "rich" enough to give "all" of $\mathcal{F}_{-\infty}^{\infty}$.) Gurevič [79] constructed a (nondegenerate) strictly stationary, finite–state, bilaterally deterministic sequence $X := (X_k, k \in \mathbb{Z})$ such that $\mathcal{T}_{\text{past}}$ and $\mathcal{T}_{\text{future}}$ are each trivial. Ornstein and Weiss [122] showed that among the strictly stationary, finite–state sequences that are "isomorphic to a Bernoulli shift," the bilaterally deterministic ones are in a certain sense "ubiquitous."

In [17], a (nondegenerate) strictly stationary, real (not discrete–state) random sequence $X := (X_k, k \in \mathbb{Z})$ is constructed which is both $\rho$–mixing (the mixing rate can be arbitrarily fast, short of $m$-dependence) and bilaterally deterministic. Burton, Denker, and Smorodinsky [54] constructed a (nondegenerate) strictly stationary, finite–state, strongly mixing ($\alpha(n) \to 0$), bilaterally deterministic sequence $X := (X_k, k \in \mathbb{Z})$. There they also posed the following question, which remains unsolved:

**Question 3**: *Does there exist a (nondegenerate) strictly stationary, finite–state, $\rho$-mixing, bilaterally deterministic sequence $X := (X_k, k \in \mathbb{Z})$?*

## 2.7. A question related to Bernoulli shifts

A theorem of Ornstein in [121] states that a given strictly stationary, finite–state random sequence is "isomorphic to a Bernoulli shift" if and only if it satisfies a certain condition of weak dependence known as the "very weak Bernoulli" condition. The terminology can be found in [152] and in other references on Ornstein isomorphism theory, and need not be given here.

Earlier, for a given strictly stationary, finite–state sequence $X := (X_k, k \in \mathbb{Z})$, Friedman and Ornstein [76] had shown that if $X$ satisfies absolute regularity (i.e. the weak Bernoulli condition — see (2.4)), then $X$ is isomorphic to a Bernoulli shift; and Smorodinsky [153] had shown that if $X$ is isomorphic to a Bernoulli shift, it need not satisfy strong mixing ($\alpha(n) \to 0$). In the 1970's D.S. Ornstein himself posed the following problem:

**Question 4**: *Suppose $X := (X_k, k \in \mathbb{Z})$ is a strictly stationary, finite–state, strongly mixing ($\alpha(n) \to 0$) random sequence; does it follow that $X$ is isomorphic to a Bernoulli shift?*

This question remains unsolved, even with the assumption of strong mixing replaced by $\rho$-mixing or even $\rho^*$-mixing. Martin [108] and Strittmatter [155] have shown that the answer is affirmative under the extra assumption of a sufficiently fast rate of convergence of $\alpha(n)$ to 0.

## 3. Markov Chains

Most of the material in this section can be found in [42, Chapter 7] and [43, Chapters 21 and 24]. Extensive further information can be found in [112], [120] and [138, Chapter 7].

### 3.1. Some basic facts

If $X := (X_k, k \in \mathbb{Z})$ is a (not necessarily stationary) Markov chain, then by the Markov property and an elementary argument, for each $n \geq 1$,

$$\alpha(n) = \sup_{j \in \mathbb{Z}} \alpha(\sigma(X_j), \sigma(X_{j+n})). \tag{3.1}$$

In the case where the Markov chain $X$ is strictly stationary, (3.1) reduces to

$$\alpha(n) = \alpha(\sigma(X_0), \sigma(X_n)). \tag{3.2}$$

Analogous comments apply to the other measures of dependence here. (In the case of $\psi'(n)$, the "sup" in (3.1) is replaced by "inf".) Such facts play a key role in the study of mixing conditions for Markov chains.

**Theorem 3.1** *Suppose $X := (X_k, k \in \mathbb{Z})$ is a strictly stationary, finite–state Markov chain. Then the following five statements are equivalent:*
*(a) $X$ is irreducible and aperiodic.*
*(b) $X$ is mixing (in the ergodic–theoretic sense).*
*(c) $\alpha(n) \to 0$ as $n \to \infty$.*
*(d) $\psi(n) \to 0$ as $n \to \infty$.*
*(e) $\rho^*(n) \to 0$ as $n \to \infty$.*

The equivalence of (a), (b), (c), and (d) is a well known, elementary consequence of the classic convergence ("equilibrium") theorem for strictly stationary, finite–state, irreducible, aperiodic Markov chains. Also, (e) $\Rightarrow$ (c) trivially, and (d) $\Rightarrow$ (e) as a special case of Theorem 3.3(7) below. The mixing rates are (at least) exponentially fast, and of course the other mixing conditions can be included here.

The next statement is a well known, elementary consequence of the classic convergence ("equilibrium") theorem for strictly stationary, countable–state, irreducible, aperiodic Markov chains.

**Theorem 3.2** *Suppose $X := (X_k, k \in \mathbb{Z})$ is a strictly stationary, countable–state Markov chain. Then the following four statements are equivalent:*
*(a) $X$ is irreducible and aperiodic.*
*(b) $X$ is mixing (in the ergodic–theoretic sense).*
*(c) $\alpha(n) \to 0$ as $n \to \infty$.*
*(d) $\beta(n) \to 0$ as $n \to \infty$.*

In the rest of section 3, the state space of the Markov chains is not necessarily countable. For convenience, the state space is taken to be $\mathbb{R}$.

**Theorem 3.3** *Suppose $X := (X_k, k \in \mathbb{Z})$ is a (not necessarily stationary) Markov chain. Then each of the following statements holds:*
   *(1) If $\rho(n) < 1$ for some $n \geq 1$, then $\rho(n) \to 0$ at least exponentially fast as $n \to \infty$.*
   *(2) If $\phi(n) < 1/2$ for some $n \geq 1$, then $\phi(n) \to 0$ at least exponentially fast as $n \to \infty$.*

*(3) If $\psi'(n) > 0$ for some $n \geq 1$, then $1 - \psi'(n) \to 0$ at least exponentially fast as $n \to \infty$.*

*(4) If $\psi(n) < 1$ for some $n \geq 1$, then $\psi(n) \to 0$ at least exponentially fast as $n \to \infty$.*

*(5) If $\psi^*(n) < 2$ for some $n \geq 1$, then $\psi(2n) < 1$ for the same $n$, and the conclusion of (4) holds.*

*(6) If $\rho^*(n) < 1$ for some $n \geq 1$, then $\rho^*(n) \to 0$ at least exponentially fast as $n \to \infty$.*

*(7) If $\psi'(n) > 0$ for some $n \geq 1$, then $X$ is $\rho^*$–mixing (and the conclusions of (3) and (6) hold).*

Statement (1) was pointed out in [138, p. 216, lines 1–3], statement (2) was pointed out in [59, Theorem 4] and is a variation on a result of Doeblin [69] (see Theorem 3.4(2) below), and statement (4) (in essence) was pointed out in [9, Lemma 8 and Theorem 5]. Statements (3) and (6) are variations on (4) and (1). Statement (5) was pointed out in [33]. Statement (7) was shown in [31, Theorem 1.2].

In connection with (7), the following question (see Question 2 in section 2.3) remains open:

**Question 5**: *If a given Markov chain is $\phi$–mixing, does it follow that it is $\rho^*$–mixing?*

The paper [36] gives an example of a strictly stationary, countable–state Markov chain which is $\rho$–mixing but not $\rho^*$–mixing. In [34, Example 6.4], [42, Example 7.16], there is an example of a strictly stationary Markov chain which is $\rho^*$–mixing but does not satisfy absolute regularity.

**Theorem 3.4** *Suppose $X := (X_k, k \in \mathbb{Z})$ is a strictly stationary Markov chain which is ergodic and aperiodic.*

*(1) If $\alpha(n) < 1/4$ for some $n \geq 1$, then $\alpha(n) \to 0$ (but not necessarily exponentially fast) as $n \to \infty$.*

*(2) If $\phi(n) < 1$ for some $n \geq 1$, then $\phi(n) \to 0$ (at least exponentially fast) as $n \to \infty$.*

Statement (1) is implicitly contained in arguments of Rosenblatt [139]. Statement (2) is a well known reformulation, in the language of strong mixing conditions, of a classic result of Doeblin [69], in connection with what is now known as "Doeblin's condition" (see section 3.2 below). A detailed exposition of both statements (1) and (2) (along with some basic details on "Doeblin's condition") is given in [43, Theorem 21.22, section 21.23, and Theorem 21.25]. (For more on "Doeblin's condition," see also [112], [120], and [138]. See also section 3.2 below.)

For (not necessarily Markovian) strictly stationary sequences, Cogburn [57] introduced the "uniform ergodicity" condition, a "Cesaro" variant of the strong mixing condition. Rosenblatt [139, Theorem 1] showed that if a strictly stationary Markov chain satisfies both uniform ergodicity and mixing (in the ergodic–theoretic sense), then it is strongly mixing. With a slight extension of Rosen-blatt's argument, it was shown in [41] that if a strictly stationary Markov chain is ergodic and aperiodic and satisfies $UERG(n) < 1/4$ for some $n \geq 1$, then it is strongly mixing. (Here $UERG(n)$, $n = 1, 2, 3, \ldots$, are the "dependence coefficients" associated with the uniform ergodicity condition.)

### 3.2. Harris recurrence, geometric ergodicity, and Doeblin's condition (again)

Suppose $X := (X_k, k \in Z)$ is a strictly stationary Markov chain. Let $\mu$ denote the (marginal) distribution of $X_0$ (on $(\mathbb{R}, \mathcal{R})$). We shall use the notation $P((X_1, X_2, X_3, \ldots) \in B | X_0 = x)$, $x \in \mathbb{R}$, $B \in \mathcal{R}^{\mathbb{N}}$, to denote a regular conditional distribution of $(X_1, X_2, X_3, \ldots)$ given $X_0 = x$.

The (strictly stationary) Markov chain $X$ is said to be "irreducible" if the following holds for $\mu$–a.e. $x \in \mathbb{R}$:

$$\forall\ B \in \mathcal{R} \text{ such that }\ \mu(B) > 0,$$
$$\exists\ n \geq 1 \text{ such that } P(X_n \in B | X_0 = x) > 0. \tag{3.3}$$

The (strictly stationary) Markov chain $X$ is said to be "Harris recurrent" [81] if the following holds for $\mu$–a.e. $x \in \mathbb{R}$:

$$\forall\ B \in \mathcal{R} \text{ such that } \mu(B) > 0,$$
$$P(X_n \in B \text{ for infinitely many } n \geq 1 | X_0 = x) = 1. \tag{3.4}$$

For broader notions of irreducibility and (Harris) recurrence, not restricted to strict stationarity, see e.g. Orey [120] or Meyn and Tweedie [112].

**Theorem 3.5** *Suppose $X := (X_k, k \in \mathbb{Z})$ is a strictly stationary Markov chain. Then the following three statements are equivalent:*
*(a) $X$ is irreducible.*
*(b) $X$ is Harris recurrent.*
*(c) $X$ is ergodic and $\lim_{n \to \infty} \beta(n) < 1$.*

*If any one (hence all three) of conditions (a), (b), and (c) hold, then for some positive integer $p$, $\lim_{n \to \infty} \beta(n) = 1 - 1/p$ and the Markov chain $X$ has period $p$ (aperiodic if $p = 1$).*

This is well known. It was shown, explicitly or implicitly, in the book by Orey [120]. As a special case, one has the following well known statement—in essence a reformulation, in the language of strong mixing conditions, of a result of Orey [119].

**Corollary 3.6** *Suppose $X := (X_k, k \in \mathbb{Z})$ is a strictly stationary Markov chain. Then the following two statements are equivalent:*
*(a) $X$ is Harris recurrent and aperiodic.*
*(b) $X$ satisfies absolute regularity ($\beta(n) \to 0$).*

In [43, Theorem 21.5 and Corollary 21.7, and Theorem 20.6], there is an exposition of Theorem 3.5 and Corollary 3.6, together with an adaptation (from Henry Berbee [5], [6], at least in spirit) of Theorem 3.5 to general (not necessarily Markovian) strictly stationary sequences.

A given strictly stationary Markov chain $X := (X_k, k \in \mathbb{Z})$ is said to satisfy "geometric ergodicity" if there exist Borel functions $a : \mathbb{R} \to (0, \infty)$ and $c : \mathbb{R} \to (0, \infty)$ such that the following holds for $\mu$-a.e. $x \in \mathbb{R}$:

$$\forall\, n \geq 1,\ \forall\, B \in \mathcal{R}, \quad |P(X_n \in B | X_0 = x) - \mu(B)| \leq a(x) \cdot e^{-c(x) \cdot n}. \qquad (3.5)$$

**Theorem 3.7** *Suppose $X := (X_k, k \in \mathbb{Z})$ is a strictly stationary Markov chain. Then the following three conditions are equivalent:*

*(a) The Markov chain $X$ satisfies geometric ergodicity.*

*(b) There exists a positive constant $c$ and a Borel function $a : \mathbb{R} \to (0, \infty)$ such that the following holds for $\mu$-a.e. $x \in \mathbb{R}$:*

$$\forall\, n \geq 1,\ \forall\, B \in \mathcal{R}, \quad |P(X_n \in B | X_0 = x) - \mu(B)| \leq a(x) \cdot e^{-cn}. \qquad (3.6)$$

*(c) The Markov chain $X$ satisfies absolute regularity with $\beta(n) \to 0$ at least exponentially fast as $n \to \infty$.*

This theorem evolved through the papers of Kendall [102], Vere–Jones [160], Nummelin and Tweedie [117], and Nummelin and Tuominen [115].

In general, in equation (3.6) in Theorem 3.7, the function $a(x)$ cannot be replaced by a positive constant. If it could, that would imply $\phi$-mixing, which for Markov chains is strictly stronger than absolute regularity with an exponential mixing rate.

For more details on geometric ergodicity, see e.g. [112]. Analogs of Theorem 3.7 involving rates of convergence slower than exponential, have been developed in numerous references, including works of Nummelin and Tuominen [116], Frenk [75], and Heinrich [82].

Now refer once more to Theorem 3.4(2). For a given strictly stationary Markov chain $X := (X_k, k \in \mathbb{Z})$, the most basic version of "Doeblin's condition" is as follows:

$$\exists\, A \in \mathcal{R} \text{ with } \mu(A) = 1,\ \exists\, \varepsilon \in (0,1),\ \exists\, n \geq 1 \text{ such that}$$
$$\forall\, x \in A,\ \forall\, B \in \mathcal{R} \text{ with } \mu(B) \leq \varepsilon, \text{ one has that}$$
$$P(X_n \in B | X_0 = x) \leq 1 - \varepsilon. \qquad (3.7)$$

By a simple argument, for the given strictly stationary Markov chain $X$, (3.7) is equivalent to the condition that for some $n \geq 1$, $\phi(n) < 1$.

For a given strictly stationary Markov chain $X$, a second version of "Doeblin's condition" is as follows: $X$ is ergodic and aperiodic and satisfies (3.7).

Theorem 3.4(2) above is an equivalent formulation, in the language of strong mixing conditions, of a classic theorem of Doeblin [69]. Doeblin's original formulation of that theorem is essentially as follows:

Suppose $X := (X_k, k \in \mathbb{Z})$ is a strictly stationary Markov chain which is ergodic and aperiodic and satisfies (3.7). Then there exists a set $A \in \mathcal{R}$ with $\mu(A) = 1$, and positive constants $Q$ and $r$, such that for every $n \geq 1$, every $x \in A$, and every $B \in \mathcal{R}$, $|P(X_n \in B | X_0 = x) - \mu(B)| \leq Qe^{-rn}$.

### 3.3. Instantaneous functions of Harris recurrent Markov chains

A (not necessarily Markovian) strictly stationary sequence $X := (X_k, k \in \mathbb{Z})$ is said to be representable as an "instantaneous function" of a strictly stationary, Harris recurrent Markov chain if the sequence $X$ has the same distribution (on $(\mathbb{R}^{\mathbb{Z}}, \mathcal{R}^{\mathbb{Z}})$) as the (strictly stationary) random sequence $(f(Y_k), k \in \mathbb{Z})$ for some (real) strictly stationary, Harris recurrent Markov chain $Y := (Y_k, k \in \mathbb{Z})$ (defined on some probability space) and some Borel function $f : \mathbb{R} \to \mathbb{R}$.

Instantaneous functions of Harris recurrent Markov chains are of interest in limit theory under dependence conditions. See e.g. the results on large deviations for such sequences in [67]. It is therefore of interest to see what strong mixing conditions (if any) might imply that kind of structure.

In [7] and [25], some (non–Markovian) strictly stationary sequences $X := (X_k, k \in \mathbb{Z})$ are constructed which are $\psi^*$-mixing (with a very fast rate of convergence of $\psi^*(n) - 1$ to 0) but which cannot be represented as a instantaneous function of a strictly stationary, Harris recurrent Markov chain. In the construction in the latter paper, the mixing rate can be made arbitrarily fast (short of $m$-dependence). As a consequence, for example, the large deviations results in [51], involving general (not necessarily Markovian) strictly stationary $\phi$-mixing sequences with a very fast mixing rate, cannot be derived as corollaries of corresponding results for instantaneous functions of strictly stationary, Harris recurrent Markov chains.

The following problem, posed in [7], remains open:

**Question 6**: *Does there exist a strictly stationary $\psi$-mixing (or perhaps even 1-dependent) sequence which cannot be represented as an instantaneous function of a strictly stationary, Harris recurrent Markov chain?*

## 4. Behavior of the dependence coefficients

We turn our attention again to general (not necessarily Markovian) strictly stationary sequences.

### 4.1. Possible limit values

For some of the dependence coefficients, there are hidden restrictions on the possible limit values. The material here in section 4.1 is treated in detail in [43, Chapter 22].

**Theorem 4.1** *Suppose $X := (X_k, k \in \mathbb{Z})$ is a strictly stationary sequence of random variables.*
  *(1) Either $\lim_{n \to \infty} \psi'(n) = 1$ ($\psi'$-mixing), or $\psi'(n) = 0 \ \forall \ n \geq 1$.*

(2) If there exists $n \geq 1$ such that $\psi^*(n) < \infty$ and $\psi'(n) > 0$, then $\lim_{n\to\infty} \psi(n) = 0$ ($\psi$-mixing).

(3) If the sequence $X$ is (strictly stationary and) mixing (in the ergodic–theoretic sense), then the following five statements hold:

(a) Either $\lim_{n\to\infty} \beta(n) = 0$ or $\beta(n) = 1 \ \forall \ n \geq 1$.
(b) Either $\lim_{n\to\infty} I(n) = 0$ or $I(n) = \infty \ \forall \ n \geq 1$.
(c) Either $\lim_{n\to\infty} \phi(n) = 0$ or $\phi(n) = 1 \ \forall \ n \geq 1$.
(d) Either $\lim_{n\to\infty} \psi^*(n) = 1$ or $\psi^*(n) = \infty \ \forall \ n \geq 1$.
(e) Either $\lim_{n\to\infty} \psi(n) = 0$ or $\lim_{n\to\infty} \psi(n) = 1$ or $\psi(n) = \infty \ \forall \ n \geq 1$.

Statements (1) and (3)(a)(b)(c)(d) are taken respectively from [14, Theorem 1 and p. 56, lines 17–26], [13, Theorem 1], [15, Lemma 0.6], [11, Theorem 1], and [14, Theorem 1]. Statement (3)(a) and its proof were a slight extension of an earlier similar statement and proof from the papers of Volkonskii and Rozanov [162], [163, Theorem 4.1, proof on pp. 194–195]. A slight adaptation of their argument was used to prove (3)(b). Special cases of (3)(c) were known earlier, for Markov chains (the result of Doeblin [69] reformulated in Theorem 3.4(2)) and for Gaussian sequences (a result of Ibragimov [89] described after Theorem 7.1). Statements (2) and (3)(e) each follow immediately from (1) and (3)(d). Statement (2) was included in [30] and [43, Theorem 22.11] at the suggestion of Manfred Denker [64].

For general (not necessarily Markovian) strictly stationary sequences, analogs of Theorem 4.1 do not hold for the dependence coefficients $\alpha(n)$, $\rho(n)$, and $\rho^*(n)$. (See Theorem 4.5 below.)

**Theorem 4.2** *Suppose $X := (X_k, k \in \mathbb{Z})$ is a strictly stationary, ergodic sequence of random variables.*

*(1) Then $\lim_{n\to\infty} \beta(n) = 1 - 1/p$ for some $p \in \{1, 2, 3, \dots\} \cup \{\infty\}$.*

*(2) Suppose the quantity $p$ in (1) satisfies $2 \leq p < \infty$. Then letting $\mathcal{J}$ denote the invariant $\sigma$–field of $T^p$ (the p-th power of the usual shift operator $T$ on events in $\mathcal{F}_{-\infty}^{\infty}$), one has the following:*

*(a) $\mathcal{J} \doteq \mathcal{T}_{\text{past}} \doteq \mathcal{T}_{\text{future}} \doteq \mathcal{T}_{\text{double}}$.*

*(b) The $\sigma$-field $\mathcal{J}$ is purely atomic, with exactly $p$ atoms, each having probability $1/p$. If $A$ is any one of those atoms, then $T^p A \doteq A$, and the $p$ atoms are $A, TA, T^2 A, \dots, T^{p-1}A$.*

*(c) Conditional on any atom of $\mathcal{J}$, the sequence of random vectors $(Y_k, k \in \mathbb{Z})$ defined by $Y_k := (X_{(k-1)p+1}, X_{(k-1)p+2}, \dots, X_{kp})$ is strictly stationary and satisfies absolute regularity.*

Theorem 4.2 is due to Henry Berbee [3, Theorem 2.1] and [6, Theorem 2.2]. (In the first of those references, an earlier preprint of [13] was cited for the related result in Theorem 4.1(3)(a).) (Here $\doteq$ means equality modulo null-sets.)

**Corollary 4.3** *Suppose $X := (X_k, k \in \mathbb{Z})$ is a strictly stationary, ergodic sequence of random variables. Suppose also that $X$ does not satisfy part (2) of Theorem 4.2 for any $p \in \{2, 3, 4 \dots\}$ — for example, suppose that for every $p \in \{2, 3, 4, \dots\}$, the invariant $\sigma$–field of $T^p$ (the p-th power of the shift operator) is trivial. Then all five statements (a)–(e) in Theorem 4.1(3) hold.*

Under the hypothesis of Corollary 4.3, statement (3)(a) in Theorem 4.1 was pointed out by Berbee [4] as an immediate consequence of Theorem 4.2; and statements (3)(c)(d)(e) in Theorem 4.1 were pointed out in [30] and [43, Corollary 22.13] (the latter also included (3)(b)) as an immediate consequence of Theorems 4.1(3) and 4.2.

Theorems 4.1 and 4.2 and Corollary 4.3 have a couple of further consequences, verified in detail in [43, Theorem 22.14 and Corollary 22.15]:

(1) Suppose $X := (X_k, k \in \mathbb{Z})$ is a strictly stationary sequence. (No assumption of ergodicity.) Then the following statements hold:

(a) If $\phi(n) < 1/2$ for some $n \geq 1$, then $X$ is $\phi$-mixing.

(b) If $\psi^*(n) < 2$ for some $n \geq 1$, then $X$ is $\psi^*$-mixing.

(c) If $\psi(n) < 1$ for some $n \geq 1$, then $X$ is $\psi$-mixing.

(2) Suppose $p \in \{2, 3, 4, \ldots\}$, and $X := (X_k, k \in \mathbb{Z})$ is a strictly stationary, ergodic sequence such that $\lim_{n \to \infty} \beta(n) = 1 - 1/p$. Then the following statements hold:

(a) Either $\lim_{n \to \infty} I(n) = \log p$ or $I(n) = \infty \ \forall \ n \geq 1$.

(b) Either $\lim_{n \to \infty} \phi(n) = 1 - 1/p$ or $\phi(n) = 1 \ \forall n \geq 1$.

(c) Either $\lim_{n \to \infty} \psi^*(n) = p$ or $\psi^*(n) = \infty \ \forall \ n \geq 1$.

(d) Either $\lim_{n \to \infty} \psi(n) = p - 1$ or $\psi(n) = \infty \ \forall \ n \geq 1$.

For random sequences that are (say) weakly stationary and strongly mixing ($\alpha(n) \to 0$) but not strictly stationary, statements above such as Theorem 4.1 do not hold in general. See e.g. the construction in [18, Theorem 7.4] (a modification of that in [103, Theorem 1]).

## 4.2. Possible mixing rates

Kesten and O'Brien [103] constructed several classes of examples that established the following general principle: For strictly stationary random sequences, the mixing rates for the various strong mixing conditions can be essentially arbitrary, and in particular, arbitrarily slow.

Later, the papers [11, Theorem 2] and [14, Theorem 2] gave respectively a couple of variations on that principle:

(1) For a given strictly stationary sequence (not $m$-dependent), the $\phi$-mixing and "time reversed $\phi$-mixing" conditions can hold simultaneously with essentially arbitrary separate mixing rates, or alternatively either condition can hold with an essentially arbitrary rate while the other fails to hold.

(2) For a given strictly stationary sequence (not $m$-dependent), the $\psi^*$-mixing and $\psi'$-mixing conditions can hold simultaneously with essentially arbitrary separate mixing rates (rates of convergence of $\psi^*(n) - 1$ and $1 - \psi'(n)$ to 0), or either condition can hold with an essentially arbitrary rate while the other fails to hold.

**Theorem 4.4** *Suppose $g : [0, \infty) \to (0, \infty)$ is a positive, continuous, strictly decreasing function such that $g(0) \leq 1/24$, $\lim_{x \to \infty} g(x) = 0$, and $\log g(x)$ is convex (as a function of $x \in [0, \infty)$). Then there exists a strictly stationary* *sequence* $X := (X_k, k \in \mathbb{Z})$ *such that for every* $n \geq 1$, $(1/4)g(n) \leq \alpha(n)$ *and* $\psi(n) \leq 8g(n)$.

That theorem was proved in [19, Theorem 1], and it was heavily based on various theorems and remarks in the paper of Kesten and O'Brien [103]. Since strong mixing ($\alpha$-mixing) and $\psi$-mixing are the "weakest" and "strongest" of the "strong mixing conditions," that theorem essentially shows that the various strong mixing conditions can occur simultaneously at essentially the same practically arbitrary rate. Thus a limit theorem under (say) $\psi$-mixing with a given mixing rate, cannot be derived from a corresponding limit theorem involving strong mixing with an essentially faster mixing rate. In a similar vein, one has the following result:

**Theorem 4.5** *Suppose* $(a_1, a_2, a_3, \dots)$, $(b_1, b_2, b_3, \dots)$, *and* $(c_1, c_2, c_3, \dots)$ *are each a nonincreasing sequence of numbers in* $[0, 1]$, *and that (i)* $4a_n \leq b_n \leq c_n$ *for all* $n \geq 1$, *and (ii)* $b_n > 0$ *for every* $n \geq 1$ *such that* $c_n > 0$. *Suppose* $(d_1, d_2, d_3, \dots)$ *is a sequence of positive numbers. Then there exists a strictly stationary sequence* $X := (X_k, k \in \mathbb{Z})$ *of random variables such that for every* $n \geq 1$,

$$a_n \leq \alpha(n) \leq a_n + d_n, \quad \rho(n) = b_n, \quad \text{and } \rho^*(n) = c_n. \tag{4.1}$$

In Theorem 4.5, the sequences $(a_n)$, $(b_n)$, and $(c_n)$ are not assumed to converge to 0.

Theorem 4.5 is a slightly modified formulation (given in [44, Theorem 26.8]) of a result in [32, Theorem 1.1]. In both references, the construction has the additional properties that (i) the marginal distribution is completely nonatomic and (ii) the double tail $\sigma$-field is trivial.

In Theorem 4.5, if (for a given $n \geq 1$) $4a_n = b_n$, then $\alpha(n) = a_n$ by (4.1) and the first inequality in (1.12). Hence for any given nonincreasing sequence $(b_1, b_2, b_3, \dots)$ of numbers in $[0, 1]$, by Theorem 4.5, there exists a strictly stationary sequence $X := (X_k, k \in \mathbb{Z})$ such that for all $n \geq 1$, $\alpha(n) = b_n/4$ and $\rho(n) = \rho^*(n) = b_n$. Corresponding general results of such "exactness" do not seem to be known for the other strong mixing conditions.

For concreteness, let us focus just on $\phi$-mixing. The following question is open:

**Question 7**: *Suppose* $(b_1, b_2, b_3, \dots)$ *is an arbitrary nonincreasing sequence of numbers in* $[0, 1]$. *Does there exist a strictly stationary random sequence* $X := (X_k, k \in \mathbb{Z})$ *such that for all* $n \geq 1$, $\phi(n) = b_n$?

In the strictly stationary $\phi$-mixing constructions in [103, Theorem 2] and [11, Theorem 2], there is an arbitrarily small "window" or "error" in the specification of the dependence coefficients $\phi(n)$. It does not even seem to be known whether there exists a strictly stationary random sequence such that $1 > \phi(1) = \phi(2) > 0 = \phi(3)$.

## 5. Independent pairs of $\sigma$-fields

Pinsker [134, p. 73] pointed out that if $(X_k, k \in \mathbb{Z})$ and $(Y_k, k \in \mathbb{Z})$ are absolutely regular sequences that are independent of each other, then the sequence $((X_k, Y_k), k \in \mathbb{Z})$ of random vectors is absolutely regular. Analogous comments apply to the other mixing conditions here. These observations are spelled out in Theorem 5.2 below, and are based on the following theorem:

**Theorem 5.1** *Suppose $\mathcal{A}_n$ and $\mathcal{B}_n$, $n = 1, 2, 3, \ldots$, are $\sigma$-fields, and the $\sigma$-fields $(\mathcal{A}_n \vee \mathcal{B}_n)$, $n = 1, 2, 3, \ldots$ are independent. Then the following statements hold:*
*(a) $\alpha(\bigvee_{n=1}^{\infty} \mathcal{A}_n, \bigvee_{n=1}^{\infty} \mathcal{B}_n) \leq \sum_{n=1}^{\infty} \alpha(\mathcal{A}_n, \mathcal{B}_n)$.*
*(b) $\rho(\bigvee_{n=1}^{\infty} \mathcal{A}_n, \bigvee_{n=1}^{\infty} \mathcal{B}_n) = \sup_{n \geq 1} \rho(\mathcal{A}_n, \mathcal{B}_n)$.*
*(c) $\beta(\bigvee_{n=1}^{\infty} \mathcal{A}_n, \bigvee_{n=1}^{\infty} \mathcal{B}_n) \leq 1 - \prod_{n=1}^{\infty}(1 - \beta(\mathcal{A}_n, \mathcal{B}_n)) \leq \sum_{n=1}^{\infty} \beta(\mathcal{A}_n, \mathcal{B}_n)$.*
*(d) $\phi(\bigvee_{n=1}^{\infty} \mathcal{A}_n, \bigvee_{n=1}^{\infty} \mathcal{B}_n) \leq 1 - \prod_{n=1}^{\infty}(1 - \phi(\mathcal{A}_n, \mathcal{B}_n)) \leq \sum_{n=1}^{\infty} \phi(\mathcal{A}_n, \mathcal{B}_n)$.*
*(e) $\psi(\bigvee_{n=1}^{\infty} \mathcal{A}_n, \bigvee_{n=1}^{\infty} \mathcal{B}_n) \leq [\prod_{n=1}^{\infty}(1 + \psi(\mathcal{A}_n, \mathcal{B}_n))] - 1$.*
*(f) $\psi^*(\bigvee_{n=1}^{\infty} \mathcal{A}_n, \bigvee_{n=1}^{\infty} \mathcal{B}_n) = \prod_{n=1}^{\infty} \psi^*(\mathcal{A}_n, \mathcal{B}_n)$.*
*(g) $\psi'(\bigvee_{n=1}^{\infty} \mathcal{A}_n, \bigvee_{n=1}^{\infty} \mathcal{B}_n) = \prod_{n=1}^{\infty} \psi'(\mathcal{A}_n, \mathcal{B}_n)$.*
*(h) $I(\bigvee_{n=1}^{\infty} \mathcal{A}_n, \bigvee_{n=1}^{\infty} \mathcal{B}_n) = \sum_{n=1}^{\infty} I(\mathcal{A}_n, \mathcal{B}_n)$.*

(Of course some of these quantities may be infinite.)

Obviously this theorem applies to the case of finitely many (say $N \geq 2$) pairs of $\sigma$-fields. (Let $\mathcal{A}_n := \mathcal{B}_n := \{\Omega, \emptyset\}$ for $n \geq N + 1$.)

Statement (b), involving the maximal correlation coefficient, is due to Csáki and Fischer [58, Theorem 6.2]. A short proof was given by Witsenhausen [165, Theorem 1]. Statement (b) is very useful in the study of the $\rho$-mixing and $\rho^*$-mixing conditions. For example, (b) played a key role in the proofs of both Theorem 3.3(7) and Theorem 4.5 as well as the constructions of most of the other $\rho$-mixing (or $\rho^*$-mixing) examples alluded to in this survey. The survey paper [34] examines the possible potential use of Theorem 5.1(b) for proving the various conjectures at the end of section 2.3 (if they are correct).

In Theorem 5.1, statements (a), (c), (d), and (h) can be found respectively in [12, Lemma 8], [16, Lemma 2.1], [11, Lemma 2.3], and [134, p. 11, Theorem 2.2(3)]; statements (f) and (g) are elementary, and (e) follows from (f) and (g).

Theorem 5.1 is also given with detailed proofs in [42, Theorems 6.1 and 6.2].

As an obvious elementary application of Theorem 5.1, one has the following theorem. Here, for a given random sequence $X := (X_k, k \in \mathbb{Z})$, the dependence coefficients $\alpha(n)$ will be denoted $\alpha(X, n)$, and analogous notations will be used for the other dependence coefficients.

**Theorem 5.2** *Suppose that for each $n = 1, 2, 3, \ldots$, $X^{(n)} := (X_k^{(n)}, k \in \mathbb{Z})$ is a (not necessarily stationary) sequence of random variables. Suppose these sequences $X^{(n)}$, $n = 1, 2, 3, \ldots$ are independent of each other. Suppose that for each $k \in \mathbb{Z}$, $h_k : \mathbb{R} \times \mathbb{R} \times \mathbb{R} \times \ldots \to \mathbb{R}$ is a Borel function. Define the sequence $X := (X_k, k \in \mathbb{Z})$ of random variables by $X_k := h_k(X_k^{(1)}, X_k^{(2)}, X_k^{(3)}, \ldots)$, $k \in \mathbb{Z}$. Then for each $m \geq 1$, the following statements hold:*
*(a) $\alpha(X, m) \leq \sum_{n=1}^{\infty} \alpha(X^{(n)}, m)$.*

(b) $\rho(X, m) \leq \sup_{n \geq 1} \rho(X^{(n)}, m)$.
(b') $\rho^*(X, m) \leq \sup_{n \geq 1} \rho^*(X^{(n)}, m)$.
(c) $\beta(X, m) \leq 1 - \prod_{n=1}^{\infty}(1 - \beta(X^{(n)}, m)) \leq \sum_{n=1}^{\infty} \beta(X^{(n)}, m)$.
(d) $\phi(X, m) \leq 1 - \prod_{n=1}^{\infty}(1 - \phi(X^{(n)}, m)) \leq \sum_{n=1}^{\infty} \phi(X^{(n)}, m)$.
(e) $\psi(X, m) \leq [\prod_{n=1}^{\infty}(1 + \psi(X^{(n)}, m))] - 1$.
(f) $\psi^*(X, m) \leq \prod_{n=1}^{\infty} \psi^*(X^{(n)}, m)$.
(g) $\psi'(X, m) \geq \prod_{n=1}^{\infty} \psi'(X^{(n)}, m)$. (Note the direction of the inequality.)
(h) $I(X, m) \leq \sum_{n=1}^{\infty} I(X^{(n)}, m)$.

If each of the functions $h_k$ is a bimeasurable isomorphism from $\mathbb{R} \times \mathbb{R} \times \mathbb{R} \times \ldots$ to $\mathbb{R}$, so that $\sigma(X_k) = \sigma(X_k^{(n)}, n \geq 1)$ for each $k \in \mathbb{Z}$, then some of the inequalities in Theorem 5.2 (the ones in (b) and (b'), and if also each $X^{(n)}$ is strictly stationary the ones in (f), (g), and (h)) will in fact be equalities.

In statement (a), if (i) for all $n \geq 1$, $\alpha(X^{(n)}, m) \to 0$ as $m \to \infty$, and (ii) for some $m \geq 1$, $\sum_{n=1}^{\infty} \alpha(X^{(n)}, m) < \infty$, then by dominated convergence, $\alpha(X, m) \to 0$ as $m \to \infty$. Similar comments (appropriately modified) apply to the other statements in Theorem 5.2.

The main applications of Theorem 5.2 involve the case where (i) each of the sequences $X^{(n)}$ is strictly stationary and (ii) the functions $h_k$, $k \in \mathbb{Z}$ are identical, and (hence, by an elementary argument) the sequence $X$ is strictly stationary. To mention just a few of many examples in the context of strict stationarity, the proof of Theorem 4.4 involved an application of Theorem 5.2(e), the proof of Theorem 4.5 involved applications of Theorem 5.2(a)(b)(b'), and the well known construction of Herrndorf [85] involved an application of Theorem 5.2(a).

Theorem 5.2 (with or without "stationarity") can be adapted trivially to a finite number, say $N$, of sequences that are independent of each other. For each $n \geq N + 1$, simply let the sequence $X^{(n)}$ be defined by $X_k^{(n)}(\omega) = 0$ for all $\omega \in \Omega$ and all $k \in \mathbb{Z}$.

## 6. Second–order properties

In this section, we shall digress and take a look at complex-valued random variables, weak stationarity, "linear dependence" conditions, and spectral density. There is a vast literature on this collection of topics, and only a tiny corner of it can be treated here.

For a further treatment of the properties of (and estimation of) spectral density under strong mixing conditions, see e.g. [101], [142], and [166].

### 6.1. CCWS random sequences and spectral density

A random sequence $X := (X_k, k \in \mathbb{Z})$ is said to be CCWS ("centered, complex, weakly stationary") if the random variables $X_k$ are complex-valued, $E|X_k|^2 < \infty$ for all $k \in \mathbb{Z}$, $EX_0 = 0$ for all $k \in \mathbb{Z}$, and there exists a function $\gamma : \mathbb{Z} \to \mathbb{C}$ such that $EX_k \overline{X}_\ell = \gamma(k - \ell)$ for all $k, \ell \in \mathbb{Z}$. Strict stationarity is not assumed.

Let $T$ denote the unit circle in the complex plane. Let $\mu$ denote normalized one-dimensional Lebesgue measure on $T$ (normalized so that $\mu(T) = 1$).

For a given CCWS random sequence $X := (X_k, k \in \mathbb{Z})$, a "spectral density function" (if one exists) is a real, nonnegative, Borel, integrable function $f : T \to [0, \infty)$ such that

$$\forall\ k \in \mathbb{Z}, \quad EX_k\overline{X}_0 = \int_{t \in T} t^k f(t) \mu(dt). \tag{6.1}$$

If one identifies each $t \in T$ with the element $\lambda \in (-\pi, \pi]$ such that $t = e^{i\lambda}$, then (6.1) takes the more familiar form

$$\forall\ k \in \mathbb{Z}, \quad EX_k\overline{X}_0 = \int_{-\pi}^{\pi} e^{ik\lambda} f(e^{i\lambda}) \frac{d\lambda}{2\pi}. \tag{6.2}$$

In the literature, the factor $1/(2\pi)$ in (6.2) is often omitted. That is not important for what follows.

For a given CCWS random sequence, a spectral density function need not exist; if it does, it is unique modulo $\mu$-null sets. (Of course every CCWS random sequence has a "spectral measure" on $T$; see e.g. [70, pp. 473-474, Theorems 3.1 and 3.2].)

### 6.2. Linear dependence coefficients

Suppose $X := (X_k, k \in \mathbb{Z})$ is a CCWS random sequence. For each positive integer $n$, define the following linear dependence coefficients:

First, let

$$r(n) := \sup \left| E \left( \sum_{k=-L}^{0} a_k X_k \right) \left( \overline{\sum_{k=n}^{n+M} a_k X_k} \right) \right| \Bigg/ \left( \left\| \sum_{k=-L}^{0} a_k X_k \right\|_2 \left\| \sum_{k=n}^{n+M} a_k X_k \right\|_2 \right) \tag{6.3}$$

where the supremum is taken over all pairs of nonnegative integers $L$ and $M$ and all choices of complex numbers $a_k$, $k \in \{-L, \ldots, 0\} \cup \{n, \ldots, n+M\}$. Here and below, $0/0$ is interpreted as 0.

Next, let

$$\zeta(n) := \sup \left| E \left( \sum_{k \in Q} X_k \right) \left( \overline{\sum_{k \in S} X_k} \right) \right| \Bigg/ \operatorname{card}(Q \cup S) \tag{6.4}$$

where the supremum is taken over all pairs of nonempty, finite, disjoint sets $Q$, $S \subset \mathbb{Z}$ such that

$$\operatorname{dist}(Q, S) := \min_{q \in Q, s \in S} |q - s| \geq n. \tag{6.5}$$

(Here and below, in (6.5), the sets $Q$ and $S$ can be "interlaced"; that is, each set can have elements between ones in the other set.)

Finally, let

$$\kappa(n) := \sup \left| E \left( \sum_{k \in Q} a_k X_k \right) \left( \overline{\sum_{k \in S} a_k X_k} \right) \right| \bigg/ \left( \left( \sum_{k \in Q} |a_k|^2 \right)^{1/2} \left( \sum_{k \in S} |a_k|^2 \right)^{1/2} \right) \tag{6.6}$$

and

$$r^*(n) := \sup \left| E \left( \sum_{k \in Q} a_k X_k \right) \left( \overline{\sum_{k \in S} a_k X_k} \right) \right| \bigg/ \left( \left\| \sum_{k \in Q} a_k X_k \right\|_2 \left\| \sum_{k \in S} a_k X_k \right\|_2 \right) \tag{6.7}$$

where in each of (6.6) and (6.7), the supremum is taken over all pairs of nonempty, finite, disjoint sets $Q, S \in \mathbb{Z}$ satisfying (6.5), and all choices of complex numbers $a_k$, $k \in Q \cup S$.

*Some observations.* Each of the dependence coefficients $r(n)$, $\zeta(n)$ $\kappa(n)$, and $r^*(n)$ is nonincreasing as $n$ increases. Also, for each $n \geq 1$,

$$r(n) \leq r^*(n) \quad \text{and} \quad \zeta(n) \leq \kappa(n). \tag{6.8}$$

(To see the second inequality, let $a_k = 1$ for each $k$ in (6.6).)

The condition $\zeta(n) \to 0$ as $n \to \infty$ is equivalent to $\kappa(n) \to 0$ as $n \to \infty$; see Theorem 6.3 below.

If $r^*(n) \to 0$ as $n \to \infty$, then $r(n) \to 0$, $\zeta(n) \to 0$, and $\kappa(n) \to 0$ as $n \to \infty$. (To see that $r^*(n) \to 0 \implies \zeta(n) \to 0$ (and $\kappa(n) \to 0$); see [40, Lemma 1.5 and Remark 1.6].)

If $\sum_{n=1}^{\infty} r(2^n) < \infty$, then $r^*(n) \to 0$ as $n \to \infty$. That implication is due to Sergey Utev [159]. A detailed exposition of it can be found e.g. in [43, Theorems 23.5–23.7].

### 6.3. Criteria for $r(n) < 1$ or $r(n) \to 0$

These criteria are provided by the following two classic theorems:

**Theorem 6.1** *For a given CCWS random sequence $X := (X_k, k \in \mathbb{Z})$ and a given positive integer $n$, the following two conditions are equivalent:*

*(a) $r(n) < 1$.*

*(b) $X$ has a spectral density function $f$ (on $T$) of the form*

$$f(t) = |p(t)|^2 \exp(u(t) + \tilde{v}(t)), \qquad t \in T$$

*where $p$ is a polynomial of degree at most $n - 1$ (constant if $n = 1$), $u$ and $v$ are real bounded Borel functions on $T$ with $\|v\|_\infty < \pi/2$, and $\tilde{v}$ is the conjugate function of $v$.*

For $n = 1$, that theorem is due to Helson and Szegö [84]. For general $n \geq 1$ it is due to Helson and Sarason [83, Theorem 6].

**Theorem 6.2** *For a given CCWS random sequence $X := (X_k, k \in \mathbb{Z})$, the following two conditions are equivalent:*
   *(a) $r(n) \to 0$ as $n \to \infty$.*
   *(b) $X$ has a spectral density function $f$ (on $T$) of the form*

$$f(t) = |p(t)|^2 \exp(u(t) + \tilde{v}(t)) \tag{6.9}$$

*where $p$ is a polynomial, $u$ and $v$ are continuous real functions on $T$, and $\tilde{v}$ is the conjugate function of $v$.*

That theorem is due to Helson and Sarason [83, Theorem 5]. (The formulation of (b) here was given by Sarason [144]). Later, Sarason [145] showed that (6.9) is equivalent to $f$ being in the class VMO ("vanishing mean oscillation").

For further perspective on those two theorems, see e.g. Peller [128] and Pourahmadi [135].

## 6.4. The mixing rate assumptions $\sum_{n=0}^{\infty} r(2^n) < \infty$ and $\sum_{n=0}^{\infty} \rho(2^n) < \infty$

For CCWS random sequences satisfying $r(n) \to 0$ as $n \to \infty$, Ibragimov [92], [93] studied in detail the connections between rates of convergence of $r(n)$ to 0 and properties of the spectral density functions. This topic is treated in detail in the book by Ibragimov and Rozanov [97, Chapter 5]. One facet of it is of particular significance to limit theory under $\rho$–mixing.

Ibragimov [90, Lemma 2], [93, Lemma 5.1] proved that if a given CCWS random sequence satisfies $\sum_{n=0}^{\infty} r(2^n) < \infty$, then it has a continuous spectral density function (on $T$). (See also [97, p. 182, Lemma 17].) That helped provide a foundation later for a central limit theorem of Ibragimov [94, Theorem 2.2] for strictly stationary, $\rho$-mixing sequences that have finite second moments and satisfy the mixing–rate assumption $\sum_{n=0}^{\infty} \rho(2^n) < \infty$. That particular mixing–rate assumption subsequently became standard in central limit theory under $\rho$-mixing with no assumption of moments of higher than second order. See for example its use in central limit theorems and weak invariance principles in [21], [123], [147], [148], and [158]. See also the almost sure invariance principle in [149] for strictly stationary, $\rho$-mixing sequences that have finite second moments and satisfy the mixing–rate assumption $\rho(n) = O((\log n)^{-(1+\varepsilon)})$ for some $\varepsilon > 0$ (just barely faster than $\sum_{n=0}^{\infty} \rho(2^n) < \infty$). Most of these results can be found in the book by Lin and Lu [107].

For the central limit theorem for strictly stationary $\rho$-mixing sequences with finite second moments, the mixing rate assumption $\sum_{n=0}^{\infty} \rho(2^n) < \infty$ is essentially as sharp as possible. That was shown with counterexamples in [20], [38]. Those constructions involved several stages, starting with stationary $\rho$-mixing Gaussian sequences with mixing rates barely slower than $\sum_{n=0}^{\infty} \rho(2^n) < \infty$. The choice of those "building block" Gaussian sequences involved delicate use of the aforementioned connections between $r(n)$ and spectral density developed by Ibragimov [92], [93] and by Ibragimov and Rozanov [97].

### 6.5. Criteria for a continuous spectral density, and related results

For CCWS random sequences, three theorems will be stated and then some related comments will be made. As indicated in the relevant sources, these theorems and subsequent comments all extend to a broader context of random fields.

**Theorem 6.3** *For a given CCWS random sequence $X := (X_k, k \in \mathbb{Z})$, the following three conditions are equivalent:*
  (a) *The sequence $X$ has a continuous spectral density function (on $T$).*
  (b) *One has $\zeta(n) \to 0$ as $n \to \infty$.*
  (c) *One has $\kappa(n) \to 0$ as $n \to \infty$.*

As stated here, this theorem was proved in [40, Theorem 1.4 and Remark 1.8]. The implication (a) $\Rightarrow$ (c) is a simple adaptation of an argument from Kolmogorov and Rozanov [104] and its extension by Rosenblatt [140], [142, pp. 73–74, Theorem 7 and Lemma 2]. The implication (c) $\Rightarrow$ (b) follows from (6.8). The implication (b) $\Rightarrow$ (a) is an adaptation of Ibragimov's [93, Lemma 5.1] proof (see also [97, p. 182, Lemma 17]) that a CCWS random sequence satisfying $\sum_{n=0}^{\infty} r(2^n) < \infty$ has a continuous spectral density (see section 6.4).

**Theorem 6.4** *For a given nondegenerate CCWS random sequence $X := (X_k, k \in \mathbb{Z})$, the following three conditions are equivalent:*
  (a) *The sequence $X$ has a (not necessarily continuous) spectral density function (on $T$) that is bounded between two positive constants.*
  (b) *One has $r^*(1) < 1$.*
  (c) *One has $r(1) < 1$ and $r^*(n) < 1$ for some $n \geq 1$.*

In that theorem, the equivalence of (a) and (b) is due to Moore [114, Theorem 1], and the equivalence of (c) with (a) and (b) is due to the author [39, Theorem 1.6].

**Theorem 6.5** *For a given nondegenerate CCWS random sequence $X := (X_k, k \in \mathbb{Z})$, the following four conditions are equivalent:*
  (a) *The sequence $X$ has a continuous positive spectral density function (on $T$).*
  (b) *One has $r^*(1) < 1$ and $r^*(n) \to 0$ as $n \to \infty$.*
  (c) *One has $r(1) < 1$ and $r^*(n) \to 0$ as $n \to \infty$.*
  (d) *One has $r(1) < 1$, $r^*(n) < 1$ for some $n \geq 1$, and $\zeta(n) \to 0$ as $n \to \infty$.*

This entire formulation was given in [39, Theorem 1.7]. Various pieces of it were contributed by the author [24, Theorem 1], [39, Theorem 1.6], [40, Theorem 1.4], the author and Utev [50, Theorem 2], Moore [114, Theorem 1], and Rosenblatt [140], [142, pp. 73–74, Theorem 7 and Lemma 2] (extending an argument of Kolmogorov and Rozanov [104]).

For a given CCWS random sequence $X := (X_k, k \in \mathbb{Z})$, the left side of (6.1) and (6.2) might be referred to as the "covariance at lag $k$." If one has a succession of CCWS random sequences that "uniformly" satisfy certain dependence conditions, and for each integer $k$ the "covariance at lag $k$" converges to a limit (a complex number), then the spectral densities may converge (either uniformly or in some weaker sense) to a limit function. In connection with Theorem 6.3, such a result was given in [40, Theorem 3.1], as a routine extension of corresponding earlier results of Falk [73] and Miller [113] involving the conditions $\sum_{n=0}^{\infty} r(2^n) < \infty$ and $r^*(n) \to 0$ respectively. In connection with Theorems 6.4 and 6.5, such results were given in Shaw [150].

We close section 6.5 with an open question:

**Question 8**: *What condition on a spectral density function is necessary and sufficient for a CCWS random sequence to satisfy $r^*(n) \to 0$ as $n \to \infty$?*

### 6.6. Interpretations involving Hilbert spaces or nonnegative definite sequences

It is well known that theorems like those here in section 6 have interpretations outside of probability theory. To cite just two references out of many, the books by Ibragimov and Rozanov [97] and Peller [128] both did much of the analysis of the linear dependence coefficients $r(n)$ in the context of Hilbert spaces not particularly tied to "probability theory." For illustrations here, let us consider the linear dependence coefficients $\zeta(n)$ and the equivalence of conditions (a) and (b) in Theorem 6.3.

(1) Suppose $H$ is a complex (or real) Hilbert space, with inner product $\langle\,.\,,\,.\,\rangle$ and norm $\|\,.\,\|$. Suppose $(h_k, k \in \mathbb{Z})$ is a sequence of elements of $H$ such that $\langle h_k, h_\ell \rangle$ depends only on $k - \ell$ ("Hilbert space stationarity"). For each $n \geq 1$, define

$$\widetilde{\zeta}(n) := \sup \left| \left\langle \sum_{k \in Q} h_k, \sum_{\ell \in S} h_\ell \right\rangle \right| \Big/ \mathrm{card}(Q \cup S)$$

where the supremum is taken over all pairs of nonempty, finite, disjoint sets $Q, S \subset \mathbb{Z}$ such that (6.5) holds. Then by applying Theorem 6.3 with an appropriate Hilbert–space isometry (or alternatively with a direct proof analogous to that of Theorem 6.3 itself), one can show that the following two conditions (a)(b) are equivalent:

(a) There exists a real nonnegative continuous function $f$ on $T$ such that for all $k \in \mathbb{Z}$, $\langle h_k, h_0 \rangle = \int_{t \in T} t^k f(t)\mu(dt)$ (where the probability measure $\mu$ on $T$ is as in (6.1)).

(b) One has $\widetilde{\zeta}(n) \to 0$ as $n \to \infty$.

(2) Suppose $(c_k, k \in \mathbb{Z})$ is a nonnegative definite sequence of complex numbers. That is, suppose that for every nonempty finite set $S \subset \mathbb{Z}$ and every choice of complex numbers $a_k$, $k \in S$ the number $\sum_{k \in S} \sum_{\ell \in S} a_k c_{k-\ell} \overline{a_\ell}$ is real and nonnegative. For each $n \geq 1$, define

$$\widetilde{\widetilde{\zeta}}(n) := \sup \left| \left\langle \sum_{k \in Q} \sum_{\ell \in S} c_{k-\ell} \right\rangle \right| \Big/ \mathrm{card}(Q \cup S).$$

where the supremum is taken over all pairs of nonempty, finite, disjoint sets $Q$, $S \subset \mathbb{Z}$ such that (6.5) holds. By Doob [70, p. 473, Theorem 3.1], there exists a CCWS random sequence $X := (X_k, k \in \mathbb{Z})$ such that $EX_k\overline{X}_0 = c_k$ for each $k \in \mathbb{Z}$. Hence by Theorem 6.3, the following two conditions (a)(b) are equivalent:

(a) There exists a real nonnegative continuous function $f$ on $T$ such that for all $k \in \mathbb{Z}$, $c_k = \int_{t \in T} t^k f(t)\mu(dt)$ (where $\mu$ is as in (6.1)).

(b) One has $\widetilde{\widetilde{\zeta}}(n) \to 0$ as $n \to \infty$.

(3) Comments (1) and (2) above can obviously be adapted to include condition (c) in Theorem 6.3. Comments (1) and (2) can obviously also be adapted to analogs of Theorems 6.1, 6.2, 6.4 and 6.5 as well as other theorems of a similar nature involving linear dependence coefficients for CCWS random sequences. With appropriate index sets (such as $\mathbb{R}$, $\mathbb{Z}^d$, or $\mathbb{R}^d$, $d \in \{1, 2, 3, \dots\}$), comments (1) and (2) can be adapted to theorems similar to Theorems 6.1—6.5 but involving CCWS random processes or random fields with index sets other than $\mathbb{Z}$ —including the original formulations of Theorems 6.3, 6.4, and 6.5 (in the relevant sources cited above), involving CCWS random fields indexed by $\mathbb{Z}^d$.

(4) Comment (1), involving Hilbert spaces, can be adapted to results on linear dependence coefficients that do not require "weak stationarity." Utev's [159] aforementioned observation that $\sum r(2^n) < \infty$ implies $r^*(n) \to 0$, involved (say) general (not necessarily weakly stationary) sequences of complex–valued, absolutely square–integrable random variables; and it has a natural analog (à la comment (1)) for complex (or real) Hilbert spaces. That was illustrated, at least in spirit, in Sergey Utev's original (unpublished) first draft of [50, Section 3], which involved abstract Hilbert spaces. With that calculation, Utev also illustrated the usefulness of the setting of an abstract Hilbert space when dealing with certain versions of linear dependence coefficients involving random variables that do not necessarily have mean 0 and may be a little awkward to deal with in the original setting of a probability space (as in the published version of [50, Section 3]).

## 7. Gaussian sequences

We return to real–valued random variables. For a stationary (real) Gaussian sequence $X := (X_k, k \in \mathbb{Z})$, the spectral density function $f$ (if it exists) is as in section 6.1, but with the term $EX_k\overline{X}_0$ in eqs. (6.1) and (6.2) replaced by $\text{Cov}(X_k, X_0)$; and (if it exists) that spectral density $f$ will satisfy $f(e^{i\lambda}) = f(e^{-i\lambda})$ for a.e. $\lambda \in [-\pi, \pi]$.

For stationary real Gaussian sequences, a thorough discussion of various mixing conditions is given by Ibragimov and Rozanov [97, Chapters 4 and 5]. Further perspective on this topic, in connection with Hankel operators, is provided by Peller [128, Chapters 8 and 9]. Here we shall just give a few basic facts.

For a given stationary mean–zero Gaussian sequence $X := (X_k, k \in \mathbb{Z})$ and a given $n \geq 1$, one has the well known equalities

$$\rho(n) = r(n) \qquad \text{and} \qquad \rho^*(n) = r^*(n) \tag{7.1}$$

by [104, Theorem 1] (combined with (6.3), (6.7), (2.2), (2.7), (1.19), and a trivial argument).

**Theorem 7.1** *Suppose $X := (X_k, k \in \mathbb{Z})$ is a nondegenerate stationary Gaussian sequence. Then the following three statements hold:*

*(1) The following three conditions are equivalent:*
*(a) $X$ is strongly mixing ($\alpha(n) \to 0$).*
*(b) $X$ is $\rho$-mixing.*
*(c) $X$ has a spectral density function $f$ (on $T$) of the form*

$$f(t) = |p(t)|^2 \exp[u(t) + \tilde{v}(t)]$$

*where $p$ is a polynomial, $u$ and $v$ are continuous real functions on $T$, and $\tilde{v}$ is the conjugate function of $v$.*

*(2) The following three conditions are equivalent:*
*(a) $X$ satisfies absolute regularity.*
*(b) $X$ satisfies information regularity.*
*(c) $X$ has a spectral density function $f$ (on $T$) of the form*

$$f(t) = |p(t)|^2 \exp \sum_{j=-\infty}^{\infty} a_j t^j$$

*(the sum converging in $\mathcal{L}^2(T)$) where $p$ is a polynomial whose roots (if there are any) lie on the unit circle and $\sum_{j=-\infty}^{\infty} |j| \cdot |a_j|^2 < \infty$.*

*(3) The following three conditions are equivalent:*
*(a) $X$ is $\phi$-mixing*
*(b) $X$ is $m$-dependent*
*(c) $X$ has a spectral density function $f$ (on $T$) of the form $f(t) = |p(t)|^2$ where $p$ is a polynomial.*

In statement (1), the equivalence of (a) and (b) is due to Kolmogorov and Rozanov [104, Theorem 2]. There they showed that for a stationary Gaussian sequence, $\rho(n) \leq (2\pi)\alpha(n)$ for all $n \geq 1$. In statement (1), the equivalence of (b) and (c) is due to Helson and Sarason [83, Theorem 5] (see also [144]) (apply (7.1) and Theorem 6.2 after the $X_k$'s are centered). In statement (2), the equivalence of (a) and (b) is due to Ibragimov and Rozanov [96], and the equivalence of (b) and (c) is due to Ibragimov and Solev [98]. Statement (3) is due to Ibragimov [89]. In fact his argument shows that for a given stationary Gaussian sequence and a given $n \geq 1$, $\phi(n) = 0$ or $1$. Obviously in statement (3), one can also list $\psi$-mixing, $\psi^*$-mixing, and $\psi'$-mixing (see remarks (a), (b), (c), (d), and (f) at the end of section 2.1).

Extending an observation of Kolmogorov and Rozanov, Rosenblatt [140], [142, pp. 73–74, Theorem 7 and Lemma 2] showed that if a stationary Gaussian sequence has a continuous positive spectral density, then it is $\rho^*$-mixing. (See (7.1) and Theorem 6.5.) By Theorem 7.1(2), one can easily construct a stationary Gaussian sequence that is $\rho^*$-mixing but not absolutely regular. Simply choose a spectral density function that is positive and continuous but sufficiently "jagged," such as $f(e^{i\lambda}) := \exp[\sum_{j=1}^{\infty} 2^{-j} \cos(4^j \lambda)]$, $\lambda \in [-\pi, \pi]$.

Refer to (7.1) and Question 8 at the end of section 6.5. The following problem remains open:

**Question 9**: *What condition on a spectral density is necessary and sufficient for a stationary Gaussian sequence to be $\rho^*$-mixing?*

## 8. A brief look at random fields

There is a large literature on strong mixing conditions for random fields indexed by $\mathbb{Z}^d$ (or $\mathbb{R}^d$) for $d \geq 2$. Here we shall briefly mention only a few recent developments on that topic. An extensive discussion on mixing conditions for random fields is given in the books by Bulinskii [53], Doukhan [71], and Zuev [168].

Various analogs of Theorem 2.1(1)(2) hold for strictly stationary random fields indexed by $\mathbb{Z}^d$ (or $\mathbb{R}^d$), as was shown in [22], [26]. As in Theorem 2.1(1)(2), those results involve both index sets being infinite in the definitions of the dependence coefficients. Those pitfalls are avoided if in the definition of the dependence coefficients, at least one of the two index sets is finite and its cardinality plays a suitable role. Indeed, in the formulation of strong mixing conditions for random fields, that has been common practice at least since the paper of Dobrushin [68]. (See e.g. [10], [60], [151] and [156].) For some examples, including ones to "separate" various "cardinality–based" strong mixing conditions, see [27], [53], [68], [71].

Another way to avoid the pitfalls analogous to Theorem 2.1(1)(2) is to impose suitable restrictions on the "shapes" of the two index sets. For example, in their adaptation of the absolute regularity condition to strictly stationary random fields indexed by $\mathbb{Z}^d$ ($d \geq 2$), Burton and Steif [55], [56] used pairs of index sets separated by an "annulus" whose inner and outer "radii" were related in a certain way.

Another recent development in limit theory for random fields under ("cardinality–based") strong mixing conditions, is the use of a "trade–off" between different "mixing rates" in the different coordinate directions of the index sets. See [105] and [49]. The paper [35] gives some examples (random fields indexed by $\mathbb{Z}^2$) to help "separate" different assumptions involving different pairs of mixing rates.

We close this survey with an open question related to "$m$-dependence" and the "tail $\sigma$-field" for random fields. In [28], for an arbitrary $d \geq 2$, a nondegenerate strictly stationary random field $X := (X_k, k \in \mathbb{Z}^d)$ is constructed such that (a) $X$ is "lattice half–space 1–dependent" and (b) modulo null–sets, $X$ is measurable with respect to its own "tail $\sigma$–field." (We refer to that paper for the technical definitions.) The random variables $X_k$ in that construction are real–valued, not discrete. After seeing a preprint of that paper, Robert Burton posed the following question, which remains open:

**Question 10**: *For a given $d \geq 2$, does there exist a nondegenerate, strictly stationary, finite–state (or even 2–state) random field $X := (X_k, k \in \mathbb{Z}^d)$ with properties (a) and (b) above?*

**Acknowledgement.** The author thanks Qi-Man Shao for his interest and encouragement on this paper, and the referee for suggestions which improved the exposition.